\documentclass[10pt]{amsart}
\usepackage{amsmath}
\usepackage{amssymb}
\usepackage{amsthm}
\input hart.sty         
\renewcommand{\>}{\rightarrow}

\begin{document}
\newtheorem{thm}{Theorem}
\newtheorem{prop}[thm]{Proposition}
\newtheorem{cor}[thm]{Corollary}
\newtheorem{defn}[thm]{Definition}
\newtheorem{lem}[thm]{Lemma}
\title{Kodaira-Iitaka Dimension on a Normal Prime Divisor}    
\author{Travis Kopp}      
\date{December 18, 2008}        
\maketitle

\section{introduction}
The Kodaira dimension of a non-singular variety is known to be an important birational invariant.
It measures the rate of growth of the number of global sections of successively larger
tensor powers of the cannonical line bundle.
Kodaira-Iitaka dimension is a generalization of Kodaira dimension to the case of a normal variety with an arbitrary line bundle. 
It is thus an invariant for measuring the positivity of a line bundle.
A Cartier divisor on a variety, for example, is called big if its corresponding line bundle has 
Kodaira-Iitaka dimension equal to the dimension of the variety.
The notion of Kodaira-Iitaka dimension can also be extended to rank one reflexive sheaves and Weil divisors as will be seen.\\
\indent
This paper is principally interested in looking at a relation between the Kodaira-Iitaka dimension of a divisor on a normal variety and the Kodaira-Iitaka dimension of related divisors on an irreducible normal subvariety of codimension $1$. 
In particular we will use geometric techniques to show a slightly broader version of the following,
\begin{thm}\label{intro main}
Let $X$ be a complete normal variety, $Y\sub X$ an irreducible complete normal Cartier divisor
and $\sL$ an invertible sheaf on $X$.
Then there exist integers $n_1>0, n_2\geq 0$ such that,
\[
\ka(X,\sL)-1\leq\ka(Y, \sL^{n_1}(-n_2 Y)\vert_Y) \]
Furthermore if $Y$ is not a fixed component of large tensor powers of $\sL$, we may take $n_1\gg n_2$.
\end{thm}
As an application we will show that Theorem \ref{intro main} implies the following inequality
when we are concerned with a subvariety of arbitrary codimension rather than a divisor only,
\begin{thm}\label{intro subvar}
Let $X$ be a complete normal variety, $A\sub X$ a complete normal subvariety
with $\codim_A(A\cap\Sing X)\geq 2$ and $\sL$ an invertible sheaf on $X$.
Then there exists a constant $\al>0$ and integers $n_1>0,\ n_2\geq0$,
such that the following inequality holds for $t\gg 0$,
\[
\al t^{\ka(X,\sL)-1}< h^0(A,\sL^{tn_1}\vert_A[\otimes]\sym^{[\, tn_2]}\sN_{A\vert X}^*) \]
(where $\sym^{[\, tn_2]}\sN_{A\vert X}^* =(\sym^{tn_2}\sN_{A\vert X}^*)^{**}$).\\
Furthermore if $A$ is not contained in the stable base locus of $\sL$, we may take
$n_1\gg n_2$.
\end{thm}
Theorem \ref{intro subvar} can be used, at least in some circumstances, to relate the Kodaira dimension of a subvariety to that of the ambient variety.
The following theorem is well known in this regard in the case of a smooth fiber of a proper fibration.

\begin{thm}[Easy addition, see e.g. \protect{\cite[\S 10]{Iitaka}}]
Let $f:X\>Y$ be a dominant proper morphism of non-singular varieties,
and let $\sL$ be an invertible sheaf on $X$.
If $Y'\sub Y$ is an open set over which $f$ is smooth, then,
\[
\ka(X,\sL)\leq\ka(X_y, \sL_y)+\dim Y\quad \text{for all}\ y\in Y', \]
\end{thm}
In \cite{PSS}
T. Peternell, M. Schneider and A.J. Sommese showed that this could be generalized
to the case of a subvariety, $A\sub X$,
not necessarily a fiber of a fibration, but with some positivity conditions on the normal bundle,
$\sN_{A\vert X}$.\\ 
\indent
For example, they produced the following inequality of Kodaira-Iitaka dimensions,

\begin{thm}[\protect{\cite[Theorem 4.1]{PSS}}] \label{PSS mt}
Let $X$ be a complete normal variety and $A\sub X$ a complete normal subvariety
with $\codim_A(A\cap\Sing X)\geq 2$.
If $\sN_{A\vert X}$ is $\Q$-effective
(i.e. $\sym^{[\, m]}\sN_{A\vert X}$ is generically generated by global sections for some $m>0$),
then,
\[
\ka(X)\leq \ka(A)+\codim_X A \]
\end{thm}
This result and other similar ones in \cite{PSS} follow from Theorem \ref{intro subvar}.\\
\indent
In \cite{PSS} the authors use an inequality, \cite[Theorem 2.1]{PSS}, similar to Theorem \ref{intro subvar} to establish their results.
However, Theorem \ref{intro subvar} has the advantage that when $A$ is not contained in the base locus of $\sL$,
the contribution of $\sN_{A\vert X}$ to the sheaves in question can be made arbitrarily small
compared to that of $\sL\vert_A$.\\
\indent
I would like to thank S\'{a}ndor Kov\'{a}cs for his advice and encouragement throughout 
the process of preparing and writing this paper.

\section{Reflexive Sheaves, Rational Maps and Kodaira Dimension}
To begin we need some preliminary observations.
In order to state our results in greater generality,
we will want to use reflexive sheaves.\\
\indent
A reflexive sheaf on a variety $X$ is a coherent sheaf $\sF$
for which the natural map, $\sF\>\sF^{**}$, to the double dual is an isomorphism.
We will be interested in the reflexive symmetric power and reflexive tensor product,
defined as follows,
\begin{defn}
For a variety $X$ and subvariety $A$ with reflexive sheaves $\sF$ and $\sG$ on $X$, let,
\[
\sym^{[\, t]}\sF =(\sym^t\sF)^{**}, \ \text{written as}\ \sL^{[\, t]}\
\text{for a rank one reflexive sheaf $\sL$}\]
\[
\text{let}\ 
\sF[\otimes]\sG =(\sF\otimes\sG)^{**} \]
\end{defn}
We will also want to use the following properties of reflexive sheaves.
Many of these are found in \cite[\S 1]{Hart2}.
\begin{prop}\ 
\begin{enumerate}
\item
The dual of any coherent sheaf is reflexive.
\item
If $X$ is a normal variety, $\sF$ a reflexive sheaf on $X$, and $U\sub X$ an open set with
$\codim_X(X\setminus U)\geq 2$, then $\sF=i_*(\sF\vert_U)$, where $i:U\>X$ is the inclusion.
Thus $H^0(X,\sF)\iso H^0(U,\sF\vert_U)$.
\item
If $X$ is a normal variety and $U$ as before with $\sL$ an invertible sheaf on $U$, then
$i_*\sL$ is a rank one reflexive sheaf on $X$.
\item
If $X$ is a non-singular variety, then any rank one reflexive sheaf on $X$ is an invertible sheaf.
\end{enumerate}
\end{prop}
We will also want some definitions to define a general class of varieties with which we'll work.
The following is inspired by the definition of small normal pairs in \cite{PSS},
\begin{defn}\ 
\begin{enumerate}
\item
A normal variety $X$ will be called {\rm almost complete} if 
for every coherent reflexive sheaf, $\sF$, on $X$,
$H^0(X,\sF)$ is finite dimensional.
In particular this will hold if
there is an embedding $X\sub\overline{X}$ of $X$ 
as an open set in a complete normal variety $\overline{X}$ with $\codim_{\overline{X}}(\overline{X}\setminus X)\geq 2$.
\item
A pair $(X,A)$ will be called a {\rm small normal pair}, if both $A$ and $X$ are almost complete normal varieties,
$A\sub X$ and $\codim_A(A\cap\Sing X)\geq 2$.
\end{enumerate}
\end{defn}
Now let $X$ be a normal variety and $U\sub X$ be its non-singular locus.
Then we have $\codim_X(X\setminus U)\geq 2$.
If $\sL$ is a rank one reflexive sheaf on $X$ then $\sL\vert_U$ is an invertible sheaf on $U$
with $H^0(X,\sL)\iso H^0(U,\sL\vert_U)$.\\
\indent
This means that given a finite dimensional linear subspace, $L$, of $H^0(X,\sL)$,
we may naturally define a linear subsystem of $\vert\sL\vert_U\vert$,
then a rational map $U\ratmap\P^N$ up to linear isomorphism of $\P^N$,
and finally a rational map $\phz_L:X\ratmap\P^N$ also up to linear isomorphism of $\P^N$.
In the case that $L=H^0(X,\sL)$ we will write $\phz_{\sL}$ for $\phz_L$.\\
\indent
Conversely if $\phz:X\ratmap\P^N$ is a rational map,
let $U$ be the intersection of the non-singluar locus of $X$ and the locus where $\phz$ is defined.
Then $\codim_X(X\setminus U)\geq 2$ and
$\phz\vert_U$ is determined by a linear subsystem of $\vert\phz^*\sO_{\P^N}(1)\vert$.
This means $\phz=\phz_L$ for some subspace $L\sub H^0(X,\phz^*\sO_{\P^N}(1))$,
where $\phz^*\sO_{\P^N}(1)=i_*(\phz\vert_U^*\sO_{\P^N}(1))$ is a rank one reflexive sheaf on $X$.
In general if $X\ratmap V$ is a rational map defined on $U\sub X$ and $\sL$ is an invertible sheaf on $V$,
we will write $\phz^*\sL$ for the rank one reflexive sheaf $i_*(\phz\vert_U^*\sL)$.  
\\
\indent
This allows us to define a notion of Kodaira-Iitaka dimension for rank one reflexive sheaves on almost complete normal varieties,
\begin{defn}[Kodaira-Iitaka Dimension]
Let $X$ be an almost complete normal variety and $\sL$ a rank one reflexive sheaf on $X$.
Define the {\rm Kodaira-Iitaka dimension}, $\ka(X,\sL)$ as follows,
\begin{enumerate}
\item
If $h^0(X,\sL^{[\, t]})=0$ for all $t>0$, let $\ka(X,\sL)=-\infty$.
\item
Otherwise, let,
\[
\ka(X,\sL) = \max\{\dim \phz_{\sL^{[\, t]}}(X)\vert t>0\}, \]
where $\phz_{\sL^{[\, t]}}(X)$ is the closure of the image of the map $\phz_{\sL^{[\, t]}}:X\ratmap\P^N$.
\end{enumerate}
\end{defn}
Kodaira-Iitaka dimension can also be expressed in terms of the assymptotic growth rate of numbers of global sections,
\begin{prop}
Let $X$ and $\sL$ be as in the definition of Kodaira-Iitaka dimension.
If $\ka(X,\sL)\neq -\infty$, then $\ka(X,\sL)$ is equal to the unique integer $\ka$
such that there exists constants $\be>\al>0$ for which,
\[
\al t^\ka < h^0(X,\sL^{[\, t]})< \be t^\ka \quad \text{or} \quad h^0(X,\sL^{[\, t]})=0 \]
for all $t>0$.
\end{prop}
\begin{proof}
This can be shown in the same way as the analagous statement for a complete normal variety and an invertible sheaf.
See e.g \cite[\S 10]{Iitaka}
\end{proof}
Note that if $L\sub H^0(X,\sL^{[\, t]})$ is a subspace for any $t>0$
then $\dim\phz_L(X)\leq\ka(X,\sL)$.\\
\indent
We will also want to consider the restriction of rational maps to prime divisors.
Let $X$ be any variety, $Y\sub X$ a prime divisor (that is, an irreducible subvariety of codimension one)
and $\phz:X\ratmap\P^N$ a rational map.
Then $\phz$ is defined on all of $X$ except for a closed set of codimension at least $2$.
In particular $\phz$ is defined on an open set of $Y$, so
there is a well defined rational map $\phz\vert_Y:Y\ratmap\P^N$ and we may speak of the closed set
$\phz(Y)\sub\P^N$ which is the closure of the image of $\phz\vert_Y$.\\
\indent
For a normal variety $X$ with non-singular locus $U\sub X$,
we may observe that the Weil divisors on $X$ restricted to $U$ are exactly the Cartier divisors on $U$,
and the rank one reflexive sheaves on $X$ restricted to $U$ are exactly the invertible sheaves on $U$.
Furthermore in each case the objects on $X$ are uniquely determined by their restriction to $U$.
Thus by the natural correspondence between Cartier divisors and invertible sheaves,
we see that there is a natural isomorphism $\Cl(X)\iso \RPic(X)$ between the divisor class group of $X$
and the group of rank one reflexive sheaves on $X$  up to isomorphism with reflexive tensor product.\\
\indent
By a slight abuse of notation we will write $\sO_X(D)$ for the reflexive sheaf corresponding to a Weil divisor, $D$,
on $X$ and $\sO_X(-D)$ for its dual.
We will say the Weil divisor $F$ is a fixed component of a rank one reflexive sheaf, $\sL$, if
$F\cap U$ is a fixed component of $\sL\vert_U$.
\begin{prop}
Let $X$ be an almost complete normal variety with rank one reflexive sheaf $\sL$.
If $\sL$ has a fixed component, $F$, then:
\[
\phz_{\sL}=\phz_{\sL(-F)} \]
\end{prop}
\begin{proof}
Since $F\cap U$ is a fixed divisor, we have $H^0(U,\sO_U(F\cap U))\iso k$.
Let $s_F\in H^0(U,\sO_U(F\cap U))$ be a generator.
Then tensor product by $s_F$ defines a map:
\[
H^0(U,\sL\vert_U(-(F\cap U)))\>H^0(U,\sL\vert_U) \]
The map must be injective and since $F\cap U$ is a fixed component of $\sL\vert_U$ it is surjective.
Since the rational maps $\phz_{\sL}$ and $\phz_{\sL(-F)}$
are defined by ratios between sections on an open set, they must be equal.
\end{proof}
Now for an almost complete normal variety $X$,
suppose $Y\sub X$ is a prime divisor which is also normal and almost complete.
Let $\sL$ be a rank one reflexive sheaf on $X$.
There is a natural restriction map:
\[
r_Y:H^0(X,\sL)\>H^0(Y,\sL\vert_Y^{**}) \] 
If $Y$ is a fixed component of $\sL$ then $\im r_Y =\{0\}$.
However if this is not the case, $\im r_Y$
defines a linear system of $H^0(Y,\sL\vert_Y^{**})$.
It is clear that the rational map corresponding to this linear system is in fact $(\phz_{\sL})\vert_Y$.
We therefore have the following:
\begin{prop}\label{kod im}
If $X$, $Y$ and $\sL$ are as above
and $Y$ is not a fixed component of $\sL$, then:
\[
\dim \phz_{\sL}(Y)\leq \dim \phz_{\sL\vert_Y^{**}}(Y)\leq \ka(Y,\sL\vert_Y^{**}) \]
\end{prop}

\section{Abstract Prime Divisors}
Next we will look more closely at prime divisors and introduce the notion of abstract prime divisors.
A prime divisor on a normal variety
determines a discrete valuation on the function field, $K(X)$.
We will abstract this to define the concept of abstract prime divisors for a variety. These are basically equivalent to algebraic valuations (see \cite{K+M}).\\
\indent  
Let $X$ be a variety over an algebraically closed field, $k$, with $\dim X= n$.
Let $Y$ be a prime divisor on $X$, such that the local ring at the generic point of $Y$, $\sO_{Y,X}$,
is integrally closed.
Then $\sO_{Y,X}$ will be a discrete valuation ring with a maximal ideal, $m_Y$,
induced by a valuation, $\nu_Y$, on $K(X)$,
the function field of $X$.
Furthermore, $\sO_{Y,X}/m_Y \iso K(Y)$.
This means $\ptd_k \sO_{Y,X}/m_Y = n-1$.
\begin{defn}
Let $X$ be variety.
If $R\sub K(X)$ is a discrete valuation ring induced by a valuation on $K(X)$
with maximal ideal $m_R$ such that:
\[
\ptd_k R/m_R = n - 1 \]
we will say $R$ is an {\rm abstract prime divisor} for $X$.\\
\indent
If $X'$ is a variety, $\phz:X'\ratmap X$ is a birational map and
$Y$ is a prime divisor on $X'$ such that $(\sO_{Y,X'},m_Y)\iso(R,m_R)$ in $K(X)\iso K(X')$
then we will say that $(X',Y)$ is a {\rm model} for $(X, R)$.
\end{defn}

It is not a priori clear that for a given abstract prime divisor, $R$, for $X$
there is a model for $(X,R)$.
However we may always associate $R$ with a closed subvariety of $X$ if $X$ is complete.\\
\indent
In this case, $X$ is proper over $\Spec k$.
Since there is natural map $\Spec K(X)\> X$,
we know by the valuative criterion of properness that $\Spec R\>\Spec k$
has a unique lift, $f:\Spec R\> X$.

\begin{defn}
We will give the name $\cen_X R$ to the closure of the image of the unique closed point in $\Spec R$
under the above morphism, $f:\Spec R\> X$.
Then $\cen_X R$ will be a closed subvariety in $X$.
\end{defn}

The following proposition is Lemma 2.45 in \cite{K+M}.
\begin{prop}\label{res ad}
Let $X$ be a complete variety over $k$,
and $R\sub K(X)$ an abstract prime divisor for $X$.
Define inductively:\\
\indent
$X_0=X$\\
\indent
If $(X_i, \cen_{X_i} R)$ is not a model for $R$, then let $\pi_{i+1}:X_{i+1}\>X_i$
be the blow-up of $X_i$ at $\cen_{X_i} R$\\
\indent
This process eventually terminates.
That is, for some $n\geq 0$ $(X_n, \cen_{X_n} R)$ is a model for $R$.
\end{prop}

\section{Restriction of an Abstract Prime Divisor to a Subfield} 
\indent
Let $k$ be an algebraically closed field.
Let $(K,R)$ be a pair, where
$K$ is a finitely generated field over $k$ with $\ptd_k K = n$
and $R\sub K$ is a discrete valuation ring with maximal ideal $m_R$ and 
$\ptd_k R/m_R = n-1$ corresponding to a valuation, $\nu$, on $K$.
The example we have in mind of course is $K=K(X)$,
the function field of a variety
with $R\sub K(X)$ an abstract prime divisor for $X$.\\
\indent
Consider a subfield, $K'\sub K$,
where $K'$ is also finitely generated over $k$ with $\ptd_k K' = m\leq n$.
It is possible that $K' \cap m_R=\{0\}$.
If this is not the case, $\nu\vert_{K'}$ is a valuation on $K'$.
Then $K'\cap R\sub K'$ is a discrete valuation ring
with maximal ideal $K'\cap m_R$.
\begin{prop}
Let $K,R,K',m$ be as above.
If $K'\cap m_R\neq\{0\}$, then, 
\[
\ptd_k (K'\cap R)/(K'\cap m_R) = m-1 \]
\end{prop}
\begin{proof}
We know $\ptd_k (K'\cap R)/(K'\cap m_R) \leq m-1$,
since $K'\cap m_R \neq \{0\}$ is a prime ideal of height at least $1$ in $K'\cap R$.\\
\indent
Let $\{f_1,\ldots, f_r\}\sub K'\cap R$ be chosen
such that $\{\overline{f_i}\}$ is a trancendence basis for 
$(K'\cap R)/(K'\cap m_R)\sub R/m_R$ over $k$.
Then complete a trancendence basis for $R/m_R$ by choosing
$\{g_1,\ldots, g_s\}\sub R$ so that $\{\overline{f_i}, \overline{g_j}\}$
is a trancendence basis for $R/m_R$ over $k$.
\begin{lem}
In the above $\{g_j\}$ are algebraically independent and trancendent over $K'$.
\end{lem}
\begin{proof}
Suppose otherwise, then there would be a nontrivial polynomial
\[
F(x_1,\ldots x_s)=\sum c_ix_1^{a_{i,1}}\cdots x_s^{a_{i,s}} \]
with coefficients in $K'$ and $F(g_1,\ldots, g_s)=0$.
Consider the valuations, $\nu(c_i)$, of the coefficients of $F$.
Let the minimum be achieved by $c_{\min}$.
Then $c_{\min}^{-1}F$ has coefficients in $K'\cap R$,
and $\overline{c_{\min}^{-1}F}$ is a nontrivial polynomial
with coefficients in $(K'\cap R)/(K'\cap m_R)$ such that
$F(\overline{g_1},\ldots,\overline{g_s})=0$.
This would make $\{\overline{g_j}\}$ algebraically dependent over $(K'\cap R)/(K'\cap m_R)$
and thus over $k(\overline{f_i})$.
But this contradicts the choice of $\{g_j\}$.
\end{proof}
The lemma implies 
\[
s\leq \ptd_{K'} K = n-m \]
Also we know 
\[
r+s=n-1=\ptd_k R/m_R \]
Thus we have $r\geq m-1$,
proving the proposition.
\end{proof}
We wish to consider the geometric consequences of this algebraic fact.
We immediately have the following.

\begin{cor}\label{rat ad}
Let $X$ and $V$ be varieties and
$\phz:X\ratmap V$ a dominant rational map
corresponding to $K(V)\sub K(X)$,
let $R\sub K(X)$ be an abstract prime divisor for $X$.
Then either $R\cap K(V)=K(V)$ or $R\cap K(V)$ is
an abstract prime divisor for $V$.
\end{cor}

Suppose in the above corollary
that $(X,Y)$ is a model for $R$.
Then $R\cap K(V)=K(V)$ implies that the generic point of $Y$ is taken to the generic point of $V$,
in other words, $\phz(Y)=V$.
Suppose however that $R\cap K(V)$ is an abstract divisor for $V$,
then we can reach the following conclusion.

\begin{prop}\label{rat im Y}
Let $X, V, \phz, R$ be as in Corollary \ref{rat ad}.
Suppose $(X,Y)$ is a model for $R$ and $V$ is complete.
If $R\cap K(V)=K(V)$ then $\phz(Y)=V$.
If $R\cap K(V)\neq K(V)$ then
$\phz(Y)=\cen_V R\cap K(V)$.
\end{prop}

\begin{proof}
The first case we have already dealt with.
In the second case we have seen that $R\cap K(V)$ is an abstract prime divisor for $V$.
Let $U\sub X$ be the domain of definition of $\phz$.
We have a map $\Spec K(X)\>U\>V$.
Then the properness of $V$ implies the following commuting diagram:
\[
\begin{array}{ccccc}
\Spec K(V)&\>&V&\leftarrow&\Spec K(X)\\
\downarrow&\nearrow&&\nwarrow&\downarrow\\
\Spec R\cap K(V)&&\longleftarrow&&\Spec R\\
\end{array}\]

The placing of the map $\Spec R\>\Spec R\cap K(V)$
is arrived at as follows:
The map $\Spec K(X)\>V$ must factor through the natural map $\Spec K(X)\> \Spec K(V)$.
Then, since the inclusion $R\cap K(V)\sub K(X)$ factors as $R\cap K(V)\sub R\sub K(X)$,
the map $\Spec K(X)\>\Spec R\cap K(V)$ must factor naturally through $\Spec R$.\\
\indent
Now $\phz(Y)\sub V$ is the closure of the image of the closed point of $\Spec R$ in $V$.
But this point is taken to the closed point of $\Spec R\cap K(V)$ in the appropriate map,
and the closure of the image of this point in $V$ is just $\cen_V R\cap K(V)$.
This gives our result.
\end{proof}

\section{Kodaira-Iitaka Dimension on a Prime Divisor}
\begin{prop}\label{main kod}
Let $X$ be an almost complete normal variety,
$Y\sub X$ a prime divisor that is also normal and almost complete.
Let $\phz:X\ratmap\P^N$ be a rational map.
Then, either $\phz(Y)=\phz(X)$,
or,
\[
\ka(Y, (\phz^*\sO_{\P^N}(n)[\otimes]\sO_X(-bY))\vert_Y^{**})\geq \dim \phz(X) -1 \quad \text{for}\ n\gg b > 0\]
\end{prop}
\begin{proof}
With the situation as above.
Let $V = \phz(X)\sub\P^N$, the closure of the image.
Then $V$ is complete
and the rational map $\phz$ correspond to the inclusion $K(V)\sub K(X)$.
Let $Y\sub X$ correspond to the abstract prime divisor $R\sub K(X)$.
To prove the theorem we may assume $\dim \phz(Y)<\dim V$,
since otherwise we are done.\\
\indent
First suppose that $\phz(Y)\sub V$ is a divisor.
Let $H$ be the divisor on $V$ corresponding to $\sO_V(1)$,
and $Z$ be a Cartier divisor on $V$ with $Z\geq \phz(Y)$.
$H$ is ample which means it has positive intersection with every class in $\overline\NE(V)$,
the cone of curves on $V$.
Since $V$ is projective, a cross section of $\overline\NE(V)$ is compact.
This means $H-\frac{b}{n}Z$ is ample as a $\Q$-divisor for $n\gg b>0$.
Thus the line bundle $(\sO_V(n)\otimes\sO_V(-bZ))\vert_{\phz(Y)}$ is ample.
It follows that,
\begin{align*}
\ka(Y,(\phz^*\sO_{\P^N}(n)[\otimes]\sO_X(-bY))\vert_Y^{**})
&\geq \ka(Y,\phz^*((\sO_V(n)\otimes\sO_V(-bZ))\vert_{\phz(Y)}))\\
&\geq \dim \phz(Y)\\
&=\dim \phz(X) -1 \end{align*}

\indent
If $\phz(Y)$ is not a divisor, then 
$R\cap K(V)\neq K(V)$.
By Corollary \ref{rat ad} and Proposition \ref{rat im Y}, 
$R\cap K(V)$ is an abstract prime divisor for $V$,
which has center $\phz(Y)$.
Then according to Proposition \ref{res ad}, there is a resolution by blow-ups,
\[
\tilde{V}=V_n\stackrel{\pi_n}{\>}V_{n-1}\>\cdots\>V_1\stackrel{\pi_1}{\>}V_0=V \]
such that $(\tilde{V}, \cen_{\tilde{V}}R\cap K(V))$ is a model for $R\cap K(V)$.
Let $W_i\sub V_i$ be the exceptional divisor of the blow-up $\pi_i$.
By a slight abuse of notation, let $W_i$ also represent the proper transform of $W_i$ in $\tilde{V}$.
Note that these $W_i$ are Cartier.

\begin{lem}\label{amp tildeV}
Let $H$ be a very ample divisor on $V$,
then we may choose $\{\al_i\in\Q\},\ 1\gg\al_1\gg\ldots\gg\al_n>0$ so that:
\[
\pi^*H-\sum\al_i W_i \]
is ample as a $\Q$-divisor on $\tilde{V}=V_n$.
\end{lem}

\begin{proof}
According to \cite[II.7.10]{Hart1} and the fact that $\sO_{V_1}(1)=\sO_{V_1}(-W_1)$,
we have that $n_1\pi_1^*H-W_1$ is very ample on $V_1$ for $n_1\gg 0$.
A slight adaptation of the proof shows that
$n_1\pi_1^*H-b_1W_1$ is also ample for $n_1\gg b_1>0$.
Thus $\pi_1^*H-\al_1 W_1$ is ample, where we can choose $\al_1\in\Q$ arbitrarily small and positive.
Continuing in this manner we get our result, where
at each step we choose $\al_i$ arbitrarily small compared to the $\al$'s already chosen.
\end{proof}

Let $\tilde{\phz}: X\>\tilde{V}$ be the appropriate rational map.
Then we know by Proposition \ref{rat im Y} that
$\tilde{\phz}(Y)=\cen_{\tilde{V}}R\cap K(V)$.
Since this is a model for an abstract prime divisor,
we know $\dim \tilde{\phz}(Y)=\dim \tilde{V}-1=\dim V -1$.\\
\indent
Since $V\sub\P^N$, $\sO_{\P^N}(1)\vert_V=\sO_V(1)$ is a very ample divisor on $V$.
Let $L$ be the divisor on $X$ corresponding to $\phz^*\sO_{\P^N}(1)$.
By Lemma \ref{amp tildeV},
$(\pi^*\sO_V(n)\otimes\sO_{\tilde{V}}(-n\sum\al_i W_i))$ is very ample on $\tilde{V}$
for $n\gg 0$ and $1\gg\al_1\gg\ldots\gg\al_n>0$.
Let $\psi:\tilde{V}\>\P^M$ be the embedding corresponding to this divisor.
Then $\psi\circ\tilde{\phz}:X\ratmap \P^M$ is a rational map
with $\dim \psi\circ\tilde{\phz}(Y)=\dim V-1$.
Note that $\tilde{\phz}^*\pi^*\sO_V(1)=\phz^*\sO_V(1)$,
since they agree on the domain of definition of $\tilde{\phz}$, $U$.
Therefore we have,
\[
(\psi\circ\tilde{\phz})^*\sO_{\P^M}(1) = \sO_X(n(L-\sum\al_i\tilde{\phz}^*W_i))\]
For each $W_i$ we have,
\[
h^0(X, \sO_X(\tilde{\phz}^*W_i))=h^0(U, \sO_U(\tilde{\phz}\vert_U^*W_i))
\geq h^0(\tilde{V}, \sO_{\tilde{V}}(W_i))>0 \]
Therefore all the $\tilde{\phz}^*W_i$ are effective.
Let $\tilde{\phz}^*W_i = D_i+m_iY$,
where $D_i$ is an effective (possibly zero) divisor on $X$ distinct from $Y$, $m_i\geq 0$ and $m_n>0$
since $\tilde{\phz}(Y)\sub W_n$.
Then $\psi\circ\tilde{\phz}$ is the rational map corresponding to a linear subsystem of,
\[
H^0\left(X, \sO_X\left(n\left(L-\sum\al_iD_i-\left(\sum \al_i m_i\right) Y\right)\right)\right) \]
Since the corresponding divisor is the pull back of $\sO_{\P^M}(1)$, it will not have any fixed components.
Therefore by Proposition \ref{kod im}, we have,
\[
\ka(Y,\sO_X\left(n\left(L-\sum\al_iD_i-\left(\sum \al_i m_i\right) Y\right)\right)\big\vert_Y^{**})\geq \dim \psi\circ\tilde{\phz}(Y)=\dim V -1 \]
Also each $D_i$ intersects $Y$ non-negatively,
so we have,
\[
\ka(Y, \sO_X\left(n\left(L-\left(\sum\al_i m_i \right) Y\right)\right)\big\vert_Y^{**})\geq \ka(Y,\sO_X\left(n\left(L-\sum\al_iD_i-\left(\sum \al_i m_i\right) Y\right)\right)\big\vert_Y^{**})\]
All of this together with the fact that the $\al_i$ can be chosen arbitrarily small,
while the $m_i$ are fixed by our resolution $\pi:\tilde{V}\>V$
gives,
\[
\ka(Y, \left(\phz^*\sO_{\P^N}(n)[\otimes]\sO_X(-bY)\right)\vert_Y^{**})\geq \dim \phz(X) -1 \quad \text{for}\ n\gg b>0 \]
\end{proof}

\begin{thm}\label{main div}
Let $X$ and $Y\sub X$ be as in Proposition \ref{main kod} and
let $\sL$ be a rank one reflexive sheaf on $X$.
If $\ka(X,\sL)\geq 0$, then for $a>0$, large and divisible, there exists an integer $m\geq 0$
depending on $a$,
which is taken to be zero if $Y$ is not a fixed component of $\sL^{[\,a]}$,
such that one of the following hold,
\[
\ka(Y, (\sL^{[\,a]}[\otimes]\sO(-Y)^{[\,m]})\vert_Y^{**})\geq \ka(X,\sL) \]
or,
\[
\ka(Y, (\sL^{[\,a]}[\otimes]\sO(-Y)^{[\,m]})\vert_Y^{**})\geq 0 \ \text{and,}\]
\[
\ka(Y,(\sL^{[\,na]}[\otimes]\sO_X(-Y)^{[\,nm+b]})\vert_Y^{**})\geq \ka(X,\sL)-1 \quad \text{for} \ n\gg b>0\]
\end{thm}

\begin{proof}
If $\ka(X,\sL)=-\infty$ we are done.
Thus we may assume that for $a$ large and divisible,
$\phz=\phz_{\sL^{[\,a]}}:X\ratmap\P^N$ is a rational map
and $\ka(X,\sL)=\dim \phz(X)$.
Let $\sL^{[\,a]}$ have maximum fixed component $F=D+mY$.
Then $\sL^{[\,a]}\otimes\sO_X(-F)$ is without fixed components
but corresponds to the same rational map as $\sL^{[\,a]}$.
We see $\sL^{[\,a]}[\otimes]\sO_X(-F)=\phz^*\sO_{\P^N}(1)$.\\

\indent
Now we apply Proposition \ref{main kod}.
We have by Proposition \ref{kod im} and the fact that $h^0(Y,\sO_X(D)\vert_Y^{**})>0$,
\begin{align*}
\ka(Y,(\sL^{[\,a]}[\otimes]\sO_X(-Y)^{[\,m]})\vert_Y^{**})
&\geq \ka(Y, (\sL^{[\,a]}[\otimes]\sO_X(-F))\vert_Y^{**})\\
&\geq\dim\phz(Y)
\end{align*}
If $\phz(Y)=\phz(X)$ we have,
\[
\ka(Y, (\sL^{[\,a]}[\otimes]\sO(-Y)^{[\,m]})\vert_Y^{**})\geq \ka(X,\sL) \]
Otherwise, we still have,
\[
\ka(Y, (\sL^{[\,a]}[\otimes]\sO(-Y)^{[\,m]})\vert_Y^{**})\geq 0 \]
Also, by Proposition \ref{main kod}, in this case,
\[
\ka(Y, \left((\sL^{[\,a]}[\otimes]\sO_X(-F))^{[\,n]}[\otimes]\sO(-bY)\right)\vert_Y^{**})\geq \dim \phz(X)-1=\ka(X,\sL)-1 \]
for $n\gg b>0$.\\
\indent
However,
\begin{align*}
\ka(Y, \left((\sL^{[\,a]}[\otimes]\sO_X(-F))^{[\,n]}[\otimes]\sO(-bY)\right)\vert_Y^{**})
&= \ka(Y,(\sL^{[\,na]}\otimes\sO_X(-nD-(nm+b)Y))\vert_Y^{**}) \\
&\leq \ka(Y,(\sL^{na}\otimes\sO_X(-Y)^{[\,nm+b]})\vert_Y^{**}) 
\end{align*}

\end{proof} 
\begin{cor}
Let $X$ and $Y\sub X$ be as in Proposition \ref{main kod}
and let $D$ be an effective big divisor on $X$
such that $Y$ is not a fixed component of $D$.
Then for $\epz\in\Q,\ 1\gg \epz > 0$, $(D-\epz Y)\vert_Y$ is big as a $\Q$-divisor on $Y$.
\end{cor}

\section{Application to Subvarieties of any Codimension}\label{apps}
By using blow-ups, we
may apply our results to the situation of any small normal pair, $(X,A)$.\\
\indent
We will start with the case when both $X$ and $A$ are non-singular and $\sL$ is an invertible sheaf.
In order to avoid an overly cumbersome statement, the result will be simplified and slightly weakened here,
but it will be easy to see how Theorem \ref{main div} applies in this situation more generally.

\begin{prop}\label{smooth subvar}
If $X$ is a non-singular almost complete variety,
$A$ a non-singular almost complete subvariety of codimension $d$
with normal sheaf $\sN_{A\vert X}$ and
$\sL$ an invertible sheaf on $X$,
then for some constant $\al>0$ and integers $n_1>0,\ n_2\geq 0$,
the following holds for $t$ large and divisible,
\[
\al t^{\ka(X,\sL)-1}<h^0(A, \sL^{tn_1}\vert_A\otimes \sym^{tn_2}\sN_{A\vert X}^*) \]
Furthermore, if $A$ is not contained in the stable base locus of $\sL$, we may take $n_1\gg n_2$.
\end{prop}
\begin{proof}
In the case $d=1$ this follows immediately from Theorem $\ref{main div}$,
taking $n_1=na,\ n_2= nm+b$ or $nm$;
thus assume $d>1$.\\
\indent
Let $\pi:\tilde{X}\>X$ be the blow-up of $X$ centered at $A$.
Let $Y\sub\tilde{X}$ be the exceptional divisor.
Then according to \cite[II.8.24]{Hart1}, $\pi:Y\>A$
is isomorphic to the projective space bundle, $\P(\sN_{A\vert X}^*)$
with $\sN_{Y\vert\tilde{X}}^*$ corresponding to $\sO_{\P(\sN_{A\vert X}^*)}(1)$.
In particular this means by \cite[II.7.11]{Hart1} that,
\[
\pi_*\sN_{Y\vert\tilde{X}}^{-k}\iso\sym^k\sN_{A\vert X}^* \]
Since $\pi:\tilde{X}\>X$ is surjective, we have,
$\ka(X,\sL)=\ka(X,\pi^*\sL)$.
Call this Kodaira-Iitaka dimension of $X$ with respect to $\sL$, $\ka$.\\
\indent By \cite[Corollary 1.7]{Hart2}
the push forward of a coherent reflexive sheaf under a proper dominant morphism of normal schemes with
equi-dimensional fibers is again a coherent reflexive sheaf.
Thus $A$ being almost complete implies that $Y$ is almost complete, and we may apply
Theorem \ref{main div} to conclude,
\[
\ka(Y,\pi^*\sL^{n_1}\vert_Y\otimes\sN_{Y\vert\tilde{X}}^{-n_2})\geq \ka(\tilde{X},\pi^*\sL)-1=\ka-1 \]
for $n_1=na$, $n_2=nm+b$ or $nm$. \\
\indent
In terms of the assymptotic growth definition of Kodaira-Iitaka dimension, this means:
\[
\al t^{\ka-1}< h^0(Y,\pi^*\sL^{tn_1}\vert_Y\otimes\sN_{Y\vert\tilde{X}}^{-tn_2}) \]
for some $\al>0$ and $t$ large and divisible.\\
\indent
By an application of the projection formula we have:
\begin{align*}
h^0(Y,\pi^*\sL^{tn_1}\vert_Y\otimes\sN_{Y\vert\tilde{X}}^{-tn_2})
&=h^0(A, \pi_*(\pi^*\sL^{tn_1}\vert_Y\otimes\sN_{Y\vert\tilde{X}}^{-tn_2}))\\
&=h^0(A,\sL^{tn_1}\vert_A\otimes\pi_*\sN_{Y\vert\tilde{X}}^{-tn_2})\\
&=h^0(A,\sL^{tn_1}\vert_A\otimes\sym^{tn_2}\sN_{A\vert X}^*)
\end{align*}
\indent
If $A$ is not in the stable base locus of $\sL$, then we may choose $a$
so that $Y$ is not a fixed component of $\pi^*\sL^a$.
Thus in this case we may choose $n_1\gg n_2 = b$ or $0$.
\end{proof}

\begin{thm}\label{main subvar}
Let $(X,A)$ be a small normal pair.
Let $\sN_{A\vert X}$ be the normal sheaf of $A$ in $X$, and
$\sL$ a rank one reflexive sheaf on $X$.
Then for some constant $\al>0$ and integers $n_1>0,\ n_2\geq 0$,
the following hold for $t$ large and divisible.
\[
\al t^{\ka(X,\sL)-1}<h^0(A, \sL^{[\,tn_1]}\vert_A^{**}[\otimes]\sym^{[\,tn_2]}\sN_{A\vert X}^*) \]
Furthermore, if $A$ is not contained in the stable base locus of $\sL$ then we may take $n_1\gg n_2$.
\end{thm}

\begin{proof}
Let $U= X\setminus (\Sing X \cup \Sing A)$
Then $U$ and $A\cap U$ are non-singular and
$\sL\vert_U$ is an invertible sheaf.
Also, since $\codim_X(X\setminus U)\geq 2$
and $\codim_A(A\setminus(A\cap U))\geq 2$,
both $U$ and $A\cap U$ will be almost complete.\\
\indent
Then we may apply Proposition $\ref{smooth subvar}$ to $A\cap U\sub U$
and $\sL\vert_U$
and observe that,
\[
h^0(A, \sL^{[\,tn_1]}\vert_A^{**}[\otimes] \sym^{[\,tn_2]}\sN_{A\vert X}^*)
=h^0(A\cap U, \sL^{tn_1}\vert_{A\cap U}\otimes \sym^{tn_2}\sN_{A\cap U\vert U}^*) \]
and,
\[
\ka(X,\sL)=\ka(U,\sL\vert_U) \]
to conclude what we desire.

\end{proof}

\section{Connections with Peternell, Schneider and Sommese}
In \cite{PSS}, T. Peternell, M. Schneider and A.J. Sommese
establish an inequality similar to the one in Theorem \ref{main subvar} \cite[Theorem 2.1]{PSS}.
They use this inequality as a basis for proving several results relating 
Kodaira-Iitaka dimension on a variety to Kodaira-Iitaka dimension on a subvariety.
These results can be proved using Theorem \ref{main subvar} in place of \cite[Theorem 2.1]{PSS}.\\
\indent
Before stating these results we need to introduce three definitions.
\begin{defn}[$\Q$-effective]
A coherent sheaf $\sF$ on an almost complete normal variety is said to be $\Q${\rm -effective}
if these exists an integer $m>0$ such that $\sym^{[\,m]}\sF$ is generically spanned by global sections.
\end{defn}
\begin{defn}[generically nef]
A coherent sheaf $\sF$ on an almost complete normal variety, $X$ is said to be {\rm generically nef}
if for some open $U\sub X$ with $\codim_X(X\setminus U)\geq 2$, $\sF\vert_U$ is locally free and nef.
\end{defn}
\begin{defn}[Kodaira dimension]
The {\rm Kodaira dimension} of an almost complete normal variety is defined to be,
\[
\ka(X)=\ka(X, \sO(K_X)) \]
\end{defn}
\begin{thm}
Let $(X,A)$ be a small normal pair
and $\sL$ a rank one reflexive sheaf on $X$. If $\sN_{A\vert X}$ is $\Q$-effective,
then,
\[
\ka(X,\sL)\leq \ka(A,\sL\vert_A^{**})+\codim_X A \]
In particular,
\[
\ka(X)\leq \ka(A)+\codim_X A \]
\end{thm}
\begin{proof}
This is an application of Theorem \ref{main subvar}.
Let $\ka=\ka(X,\sL)$ and $d=\codim_X A$.
For $t$ large and divisible we have,
\[
\al t^{\ka-1}\leq h^0(A,\sL^{[\,tn_1]}\vert_A^{**}[\otimes]\sym^{[\,tn_2]}\sN_{A\vert X}^*)\]
for some $\al>0$, $n_1>0$ and $n_2\geq 0$.\\
\indent
If $n_2=0$, we have,
\[
\al t^{\ka-1}<h^0(A,\sL^{[\,tn_1]}\vert_A^{**}) \]
This would mean $\ka-1=\ka(X,\sL)-1\leq\ka(A,\sL\vert_A^{**})$,
which establishes the theorem.\\
\indent
Thus assume $n_2>0$.
Since $\sN_{A\vert X}$ is $\Q$-effective, 
we may assume that for $t$ large and divisible enough,
$\sym^{[\,tn_2]}\sN_{A\vert X}$ is generically spanned by global sections.
This means $\sym^{[\,tn_2]}\sN_{A\vert X}^*\monoto\sO_A^{\oplus M}$,
where $M=\rk\sym^{[\,tn_2]}\sN_{A\vert X}<\be_1(tn_2)^{d-1}$.
Then we have,
\begin{align*}
h^0(A,\sL^{[\,tn_1]}\vert_A^{**}[\otimes]\sym^{[\,tn_2]}\sN_{A\vert X}^*)
&\leq h^0(A,\sL^{[\,tn_1]}\vert_A^{**}[\otimes]\sO_A^{\oplus M})\\
&= Mh^0(A,\sL^{[\,tn_1]}\vert_A^{**})\\
&< \be_1(tn_2)^{d-1}\be_2(tn_1)^{\ka(A,\sL\vert_A^{**})}
\end{align*}
Putting this together we have,
\[
\al t^{\ka-1}<\be_1\be_2n_1^{\ka(A,\sL\vert_A^{**})}n_2^{d-1}t^{\ka(A,\sL\vert_A^{**})+d -1} \]
for $t$ large and divisible.\\
\indent
This means, 
\[
\ka \leq \ka(A,\sL\vert_A^{**})+d \]
The second statement in the theorem follows since for $\sN_{A\vert X}$ $\Q$-effective
we have, 
\[
\ka(A,\sO(K_X)\vert_A^{**})\leq \ka(A) \] 
\end{proof}
The next two results follow from Theorem \ref{main subvar}, \cite[Lemma 3.8]{PSS}, \cite[Lemma 3.9]{PSS} and \cite[Lemma 3.10]{PSS}.
\begin{thm}
Let $(X,A)$ be a small normal pair.
Suppose $A$ lies as an open subset in a complete normal variety $\overline{A}$
with $\codim_{\overline{A}}(\overline{A}\setminus A)\geq 2$, such that
$\overline{A}$ has at worst terminal singularities and admits a good minimal model.
If $\sN_{A\vert X}$ is generically nef, then,
\[
\ka(X)\leq \ka(A)+\codim_X A \]
\end{thm}
\begin{proof}
See \cite[Theorem 4.3]{PSS}.
Theorem \ref{main subvar} may be used in place of \cite[Theorem 2.1]{PSS} in the proof.
\end{proof}
\begin{thm}
Let $(X,A)$ be a small normal pair and $\sL$ a rank one reflexive sheaf on $X$.
Suppose $\sL$ is $\Q$-Cartier and $\sL\vert_A^{**}$ is semi-ample.
If $\sN_{A\vert X}$ is generically nef, then
\[
\ka(\sL)\leq \ka(\sL\vert_A^{**})+\codim_X A \]
\end{thm}
\begin{proof}
See \cite[Theorem 4.4]{PSS}.
Theorem \ref{main subvar} may be used in place of \cite[Theorem 2.1]{PSS} in the proof.
\end{proof}
To give more detail we may look at the inequality of interest in \cite{PSS}.
\begin{thm}[\protect{\cite[Theorem 2.1]{PSS}}]\label{PSS ineq}
Let $(X, A)$ be a small normal pair and $\sL$ a rank one reflexive sheaf on $X$.
Then there is a positive integer $c$ such that for all $t\geq 0$,
\[
h^0(X,\sL^{[\,t]})\leq \sum_{k=0}^{ct}h^0(A,\sL^{[\,t]}\vert_A^{**}[\otimes]\sym^{[\,k]}\sN_{A\vert X}^*)\]
\end{thm}
A similar inequality can be shown as a corollary of Theorem \ref{main div} using the methods in section \ref{apps}.
The following inequality has the advantage that the contribution from $\sN_{A\vert X}^*$
can be taken to be arbitrarily small compared to that of $\sL\vert_A^{**}$
when $A$ is not in the stable base locus of $\sL$.
\begin{cor}[of Theorem \ref{main div}]
Let $(X, A)$ be a small normal pair and $\sL$ a rank one reflexive sheaf on $X$.
Then there exists a constant $\al>0$ and positive rational number $c$, which can be chosen to be arbitrarily
small if $A$ is not in the stable base locus of $\sL$,
such that for all large and divisible $t$,
\[
h^0(X,\sL^{[\,t]})\leq 
\al \sum_{k=0}^{\flr{ct}}h^0(A,\sL^{[\,t]}\vert_A^{**}[\otimes]\sym^{[\,k]}\sN_{A\vert X}^*)\]
\end{cor}
\begin{proof}
This follows by using the full extent of Theorem \ref{main div},
applying the blowing-up methods of section \ref{apps}
and doing some arithmetic and calculus to fit the results into the form of the sum in the corollary.
\end{proof}

\bibliographystyle{plain}
\bibliography{KodBib}

\end{document}